# An Iterative Cyclic Algorithm for Designing Vaccine Distribution Networks in Low and Middle-Income Countries

Yang Y[1], Rajgopal J[2]

**Abstract** The World Health Organization – Expanded Programme on Immunization (WHO-EPI) was developed to ensure that all children have access to common childhood vaccinations. Unfortunately, because of inefficient distribution networks and cost constraints, millions of children in many low and middle-income countries still go without being vaccinated. In this paper, we formulate a mathematical programming model for the design of a typical WHO-EPI network with the goal of minimizing costs while providing the opportunity for universal coverage. Since it is only possible to solve small versions of the model optimally, we describe an iterative heuristic that cycles between solving restrictions of the original problem and show that it can find very good solutions in reasonable time for larger problems that are not directly solvable.

## 1 Introduction

Infectious diseases are one of the main causes of mortality in low and middle-income countries, and over the years it has been repeatedly proven that one of the best strategies to protect against many diseases is an effective program of childhood vaccination. Over the years many diseases such as polio and small pox have been virtually eliminated and the incidences of many others like measles, rubella, tetanus, pertussis, and diphtheria have been vastly reduced. A World Health Organization (WHO) report on immunization (WHO, 2018a) indicates that immunization is one of the most cost-effective investments and results in 2 to 3 million deaths being averted each year.

The first major initiative to provide universal access to all important vaccines for children was the Expanded Programme on Immunization (EPI), which was established in 1974 (Bland and Clements, 1997). The second major initiative in this regard was the Global Alliance for Vaccines and Immunization, or Gavi as it is better known, which was formed in 2000 (Gavi, 2019). The main focus of Gavi was on the poorest countries, mainly in sub-Saharan Africa and Asia, where the vaccination rates were the lowest. It is estimated that EPI and Gavi have been responsible for saving millions of lives worldwide and for almost completely eliminating diseases like polio and measles (WHO, 2013).

While immunization rates have been rising steadily over the years, coverage is still only around 85% (and even lower in many countries) and it is estimated that another 1.5 million deaths per year could be avoided if global vaccination coverage could improve further. For example, in the most recent year for which data is available almost 20 million children under the age of one worldwide did not receive the recommended three doses of DTP-3 vaccine that protects against diphtheria-tetanus-pertussis (WHO, 2018b). Coverage rates are generally quite high in developed countries, but many low and middle-income countries still suffer from inadequate coverage and present vast opportunities for improvement (Gavi, 2017). There are many reasons why vaccination rates are low in poorer countries. These include limited

[1]Yuwen Yang (e-mail: yuy49@pitt.edu)
[2]Jayant Rajgopal (✉ e-mail: rajgopal@pitt.edu)
Dept. of Industrial Engineering, University of Pittsburgh, Pittsburgh, PA 15261, USA.

International Joint Conference on Industrial Engineering and Operations Management- ABEPRO-ADINGOR-IISE-AIM-ASEM (IJCIEOM 2019). Novi Sad, Serbia, July 15-17th health budgets and resources, weak or unreliable infrastructure, poorly designed distribution chains, a lack of scientific healthcare management, inadequate or faulty equipment, lack of transportation resources, poor monitoring and supervision, inadequate access to facilities, and even religious reasons (Yadav, et al., 2014). In general, the primary challenge is not so much in obtaining the vaccines as it is in ensuring that these are shipped, stored and delivered to the recipients at the end of the distribution chain in a cost-efficient fashion. A major consideration for the distribution system is that vaccines require narrowly defined temperatures of between 2 and 8°C during storage and transportation, i.e., we must deal with a so-called "cold" chain.

In most low and middle-income countries, vaccines are usually distributed via a hierarchical legacy medical network determined by political boundaries and history. These networks also have a rigid structure that is replicated in most of these countries. In particular, there has been virtually no attempt to design vaccine distribution networks with flexible structures and operational characteristics tailored to the demographic, geographical and other features of a specific country. The goal of the work presented in this paper is to introduce a mathematical model for redesigning the EPI vaccine distribution network in any country. The model adheres to established WHO guidelines and also emphasizes the operational simplicity that is important in low and middle-income countries that often lack management abilities required for more sophisticated systems. We illustrate how our model can be used by using data that is derived from several different countries.

## 2  Problem Development

A typical vaccines distribution network uses a four-tier hierarchical architecture such as the one shown in Fig. 1. Vaccines are purchased or obtained through donations from international organizations and shipped by air once or twice a year. They are then stored in a *national* distribution center in a large city (usually the capital) in the nation where it is to be distributed. Quarterly shipments of vaccines in the required amounts are then made to *regional* distribution centers (typically, fewer than 10), usually via specialized vehicles like temperature-controlled cold trucks. The regional centers then transport vaccines to *district* centers (typically, 40 to 80) that fall under their domain. This is usually done every month using 4×4 trucks with cold storage boxes or vaccine containers. The last node in the distribution chain is a *clinic* or health center (typically, several hundred of these) where infants and children come and receive vaccinations. A clinic obtains its required amounts of vaccines from the district center to which it is assigned. This usually happens once a month using locally available means of transportation such as cars, trucks, motorcycles, bicycles, etc. that move the vaccines in a vaccine carrier/cooler for storage in refrigerators at the clinic.

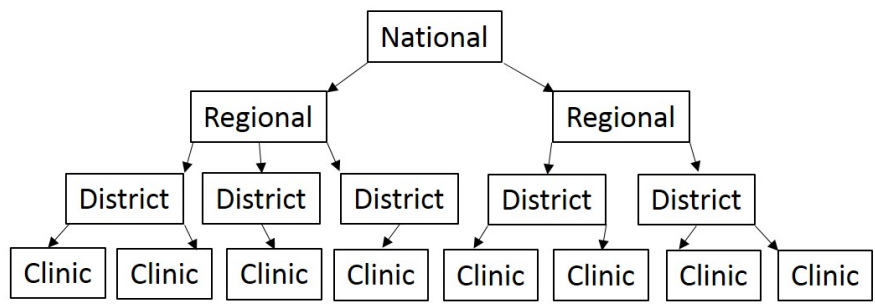

**Fig. 1.** A typical tour-tier vaccine distribution network

While there might be one more or one fewer intermediate level in some places, almost all countries distribute vaccines from the central store to individual clinics using a network that is very similar to the one just described. However, there is no compelling reason as to why the network should have this strict hierarchical structure everywhere. For instance, it might make more sense for a district to be supplied by another district, or for a clinic to be supplied directly by a regional center or for a district to be supplied by the national center. Replenishment frequencies at districts and regions could also be more flexible as long as they are not too complicated for the logistics expertise at hand.

To develop an optimal network design with more flexible operational policies, we combine the district and regional centers into the same pool and simply refer to them as hubs. Thus, our overall problem is to design a network that will move vaccines from the central store to individual clinics around the country through a system of intermediate hubs, while minimizing the overall cost of transportation, facilities and storage across the network. In line with WHO guidelines, the constraints of this optimization problem must ensure that there is *opportunity* for 100% of all demand for the area served by a clinic to be satisfied, while using recommended replenishment frequencies, picking from approved storage devices (WHO, 2018c) and available transportation modes, and meeting constraints on both storage and transportation capacity. The decisions consist of hub allocations for each clinic, all hub-to-hub and national store-to-hub connections, vaccine flows along these connections, and the types of storage and transportation devices to be used at each selected location and along each selected connection. In our problem formulation we only consider the total volume of vaccines shipped; these volumes are determined by the space requirements for the individual vaccines in the country's vaccine regimen as determined from available data. Each open facility has an associated annual operational cost. With a view to operational simplicity, which is critical in low and middle-income countries, we assume that each facility is supplied by a single node, and uses a single type of approved storage device. Each hub is supplied either by the national center or by another hub and has a flexible vaccine stock replenishment policy with either a monthly or a quarterly replenishment frequency (as per WHO guidelines). Clinics are supplied every month by a hub (or possibly by the national center if it is closer than any hub). Finally, we assume that there is a 25% safety buffer at each location as per WHO guidelines so that the total required storage volume is inflated by a factor of 1.25.

Variant of this network design problems have been addressed by the operations research community, as the *p*-median problem and the facility location problem (see Melkote and Daskin, 2001; Şahin and Süral, 2007; Melo et al., 2009). However, none of the prior work considers multiple, *discrete* choices of storage and transportation devices and replenishment frequencies, as is the case here. As far as vaccine distribution networks, Chen et al. (2014) were the first to develop a planning model. The goal was to maximize coverage under current network capacity with an extension that allowed capacity expansion. However, this work does not address the network design; rather it looks at operations in an existing network. Some relatively recent work by one of the authors of this paper and his colleagues (Lim et al., 2017) was the first to introduce a model to completely redesign the vaccine distribution network. The model we present below is based on this model but refines it considerably by dispensing with several assumptions made by the former such as coordinated deliveries to hubs using vehicle loops, fixed storage devices, rigid replenishment schedules and the option of multiple vehicle trips for a single delivery. We introduce this model next.

## *2.1 Model Formulation*

Consider a network with possible *h*ub locations in index set *H*, and *c*linics locations in index set *C*; let the index 0 to refer to the national store. Let the transportation *m*odes and WHO-qualified storage *d*evices available be indexed by index sets *M* and *D*, respectively. We use *f*=1 to denote a monthly replenishment policy with $n_1$=12 annual replenishments and *f*=2 to denote a quarterly replenishment policy with $n_2$=4 annual replenishments. Finally, define the set *A* of feasible arcs in the network, i.e., [*i-j*]∈*A* if a direct flow is possible from location *i* to location *j*. Our inputs include the following parameters and costs:

$a_i$ = volume associated with the annual demand at clinic *i*∈*C*

$V_m$ = maximum volume that can be stored in vehicle associated with transportation mode *m*∈*M*

$S_d$ = maximum volume that can be stored in storage device *d*∈*D*

$T_{ijm}$ = Roundtrip transportation cost from location *i* to location *j* using transportation mode *m*

$F_{id}$ = Facility cost per year at location *i* when it uses storage device *d*

The primary decisions that we need to make are (a) which hub locations to open, along with the storage device and replenishment frequency to use at each open location, and (b) which location $i$ to use for supplying vaccine to any location $j$ that represents a clinic or an open hub, along with the associated transportation mode. We therefore define the following decision variables and formulate the mixed integer programming model that follows

$Y_{ijmf}$ =1 if vaccines flow along [$i$-$j$] using transport mode $m$ and replenishment frequency $f$; 0 otherwise
$Z_{idf}$ = 1 if hub location $i \in H$ is open with device $d$ and replenishment frequency $f$; 0 otherwise
$X_{ij}$ = Volume of the annual vaccine flow along arc [$i$-$j$]

*Program 1*

$$\text{Minimize } \sum_i \sum_d F_{id}\left(\sum_f Z_{idf}\right) + \sum_{[i-j] \in A} \sum_m T_{ijm}\left(\sum_f n_f Y_{ijmf}\right) \quad (1)$$

s.t.

$$\sum_{i \in \{0\} \cup H} \sum_m Y_{ijm1} = 1 \ \forall j \in C \quad (2)$$

$$\sum_{i \in \{0\} \cup H} \sum_m \sum_{f \in F} Y_{ijmf} \leq 1 \ \forall j \in H \quad (3)$$

$$\sum_d \sum_f Z_{idf} \leq 1 \ \forall i \in H \quad (4)$$

$$\sum_d Z_{jdf} - \sum_{i \in \{0\} \cup H} \sum_m Y_{ijmf} = 0 \ \forall j \in H, f = 1,2 \quad (5)$$

$$\sum_{i \in \{0\} \cup H} X_{ij} = a_j \ \forall j \in C \quad (6)$$

$$\sum_{i \in \{0\} \cup H} X_{ij} - \sum_{k \in H \cup C} X_{jk} = 0 \ \forall j \in H \quad (7)$$

$$\sum_m V_m\left(\sum_f n_f Y_{ijmf}\right) \geq X_{ij} \ \forall [i-j] \in A \quad (8)$$

$$\sum_d S_d\left(\sum_f n_f Z_{jdf}\right) \geq 1.25 \sum_{i \in \{0\} \cup H} X_{ij} \ \forall j \in H \quad (9)$$

$$X_{ij} \geq 0, Y_{ijmf} \in \{0,1\}, Z_{idf} \in \{0,1\}$$

The objective (1) minimizes the sum of the annual facility costs across all facilities that are open and annual round-trip transportation costs across all arcs in the network that are selected. Constraints (2) ensure that every clinic is replenished monthly by a single hub (or the national store) with an associated transport device, while constraints (3) ensure that any open hub is supplied by a single source with an associated transport device and replenishment frequency. Constraints (4) ensure that any open hub has a single type of storage device and uses a single replenishment frequency. Constraints (5) ensure that if a hub is not open there is no flow into it, constraints (6) state that the total annual inflow into a clinic is equal to its annual demand and constraints (7) are flow balance equations for the hubs. Constraints (8) ensure that in any arc along which there is a flow, a transportation mode with sufficient capacity for that flow is selected. Lastly, constraints (9) ensure that in each open facility there is storage device with enough capacity to store the vaccines required within each replenishment interval (including the 25% buffer that the WHO requires). Note that $n_f$ is either 4 or 12 and in (8) and (9) the RHS is effectively divided by $n_f$ for the selected frequency.

### 2.1.1 Computational Limits with Program 1

We first explored the solution of Program 1 for a number of problems using a standard off-the-shelf mixed integer programming solver (IBM ILOG CPLEX 12.6) on a 3.20 GHz processor with 8 GB of memory. Problems of various sizes were generated using data that were derived from information we could access for four different countries in sub-Saharan Africa. While all of these countries currently have a similar four-tiered distribution architecture, they had significantly different demographic characteristics (size, population density, etc.) and also differed in the number of potential hub locations. In general, the effort required for a problem depends mainly on the total number of nodes as well as potential hub locations; however, it also depends on the population distribution and the costs. There was no clear means to specify the limits to what is solvable. Our numerical tests are described in more detail in Section 4, but as a gross generalization, only problems with under about 200-250 nodes and 15-20 potential hub locations can be solved in reasonable time. A key fact that makes Program 1 hard to solve is that it has a large number of 0-1 decision variables. For a typical small to mid-size 100-node problem with 15 candidate hubs, the number of binary decision variables is close to 10,000. This clearly calls for heuristics or other approaches.

## 3  An Iterative Cyclic Algorithm

We describe an easy-to-implement MIP-based heuristic that solves a sequence of MIP problems, each of which is a restricted version of Program 1 that is relatively easy to solve. These restrictions are with respect to either the replenishment frequencies used at hubs or the total number of hubs that are open. The method is motivated by initial experiments where we tested the MIP when some variables are fixed, thereby reducing the number of decision variables. First, we formulated a restricted version by fixing *all* replenishments at hubs to be done either once a month or once a quarter. These restricted versions yielded solutions in a very short amount of time and with values under 1% larger the true optimum for smaller problems. Next, we formulated an alternative restriction where we fixed the total number of hubs to be open. Once again, fixing a portion of the network structure generally yielded solutions much more quickly (although as we force more hubs to be open the time does start to increase).

Based on these observations, our algorithm starts by solving a restricted version of Program 1 with an initial vector of fixed replenishment frequencies at the hubs, to obtain a locally optimal set of open locations under this frequency vector. The algorithm then fixes these open hub locations and solves another restricted version of Program 1 (with other hubs kept closed) to find the corresponding optimal frequencies. The procedure iterates until we cannot improve the solution. Before describing the algorithm, let us denote:

$\varepsilon$: a suitably small constant

$\boldsymbol{f}$: a vector of order $|H|$ indicating the replenishment frequency at hubs; if the $i$th element is 1 then hub $i$ is set to be replenished quarterly and a constraint $\sum_{d \in D} Z_{id2} = 0$ is added to the model; alternatively if the $i$th element is 2 then hub $i$ is set to be replenished monthly, and a constraint $\sum_{d \in D} Z_{id1} = 0$ is added; if a hub is closed the corresponding element is set to 0.

$\boldsymbol{l}$: a binary vector of order $|H|$ indicating the status of each hub; if the $i$th element equals 0, hub $i$ is forced to be closed and we add a constraint $\sum_{d \in D} \sum_{f \in F} Z_{idf} = 0$ to the model; if the $i$th element is equal to 1, hub $i$ is set to be open, and we add a constraint $\sum_{d \in D} \sum_{f \in F} Z_{idf} = 1$ to the model.

$W_f^k$: locally optimum objective value at step $k$ when frequencies are fixed

$W_l^k$: locally optimum objective value at step $k$ when locations are fixed

**STEP 1**: *Initialization*

Generate a random vector of length $|H|$ where every entry is on of either 1 or 2 and define it to be $\boldsymbol{f}^1$. Note that initially, every hub is allowed to be open, and if it is open it must use the replenishment frequency specified via $\boldsymbol{f}^1$. Let $k = 1$.

**STEP 2**: *Local optimum with fixed frequencies*

Set $\boldsymbol{f} \leftarrow \boldsymbol{f}^k$ and solve Program 1 under this fixed frequency vector with the corresponding additional constraints. Let $W_f^k$ be the local optimum value obtained, with corresponding hub locations defined by the vector $\boldsymbol{l}^k$. If $k = 1$ go to STEP 3 after deleting all the constraints added at this step. Otherwise, if $W_l^{k-1} - W_f^k \leq \varepsilon$, i.e., there is no improvement, stop the algorithm with objective value $W_f^k$. Otherwise, it means that the algorithm is still improving the solution, so we delete the additional constraints added in this step and continue on to STEP 3.

**STEP 3**: *Local optimum with fixed open locations*

Set $l \leftarrow l^k$ and solve Program 1 under this fixed location vector with the corresponding additional constraints, let $W_l^k$ be the local optimum value obtained, with corresponding replenishment frequencies defined by the vector $f^k$. In the case that a hub (say, the $i^{th}$) is not open in the solution obtained at STEP 2, the $i^{th}$ element of $f^k$ is set to 0. If $W_f^k - W_l^k \leq \varepsilon$, there is no improvement; stop the algorithm with value $W_l^k$. Otherwise, delete the additional constraints added in this step and go to STEP 4.

**STEP 4:** *Update frequency*

Set $k \leftarrow k + 1$ and update the frequency vector via $f^{k+1} \leftarrow f^k$. Then return to STEP 2.

## 4 Computational Results

Table 1 below summarizes computational results from experiments that we conducted.

**Table 1** Computational Results

| No. | Hubs | Nodes | 0/1 Variables | Pop. Density | CPU: Optimum | CPU: Algorithm | Gap % |
|---|---|---|---|---|---|---|---|
| 1 | 1 | 10 | 68 | sparse | <1s | <1s | 0% |
| 2 | 1 | 18 | 116 | moderate | <1s | <1s | 0% |
| 3 | 2 | 11 | 148 | sparse | <1s | <1s | 0% |
| 4 | 3 | 22 | 420 | sparse | 2s | <1s | 0% |
| 5 | 2 | 49 | 604 | sparse | <1s | <1s | 0% |
| 6 | 3 | 39 | 726 | moderate | 2s | <1s | 0% |
| 7 | 4 | 44 | 1,088 | sparse | 3s | <1s | 0% |
| 8 | 4 | 48 | 1,184 | moderate | <1s | <1s | 0% |
| 9 | 4 | 55 | 1,352 | moderate | 4s | <1s | 0% |
| 10 | 4 | 64 | 1,568 | moderate | 7s | <1s | 0% |
| 11 | 4 | 65 | 1,592 | dense | 8s | <1s | 0.27% |
| 12 | 5 | 77 | 2,350 | moderate | 1.6s | <1s | 0% |
| 13 | 8 | 56 | 2,752 | sparse | 16s | 4s | 0% |
| 14 | 7 | 99 | 4,214 | dense | 4.4s | <1s | 0% |
| 15 | 11 | 96 | 6,424 | moderate | 10s | 5s | 0% |
| 16 | 10 | 117 | 7,100 | moderate | 146s | 29s | 0.56% |
| 17 | 14 | 101 | 8,596 | sparse | 119s | 32s | 0.28% |
| 18 | 12 | 128 | 9,312 | dense | 116s | 29s | 0.17% |
| 19 | 8 | 206 | 9,952 | dense | ~10h | 43s | 0.17% |
| 20 | 14 | 148 | 12,544 | moderate | 103s | 49s | 0% |
| 21 | 17 | 141 | 14,518 | moderate | 79s | 50s | 0% |
| 22 | 16 | 162 | 15,680 | dense | 1,304s | 46s | 0.46% |
| 23 | 13 | 210 | 16,484 | moderate | ~1d | 76s | 0.52% |
| 24 | 14 | 235 | 19,852 | dense | ~2d | 83s | 0.62% |
| 25 | 19 | 176 | 20,216 | moderate | 4,649s | 207s | 0.24% |
| 26 | 20 | 295 | 35,560 | moderate | 387s | 89s | 0.32% |
| 27 | 26 | 333 | 52,156 | moderate | 2,748s | 126s | 0.54% |

We tested our algorithm using a number of problems; as stated in Section 2.1.1 these were generated from the data we had for four different LMICs. Table 1 lists results for the 27 test problems that we were able to solve optimally. For each problem we list the number of potential hub locations, the total number of nodes in the network, the number of binary variables in the formulation of Program 1, and a label that identifies the population in the area as being dense, moderate or sparse. The number of nodes and potential hub locations in these problems ranged from 10 to 333, and from 1 to 26, respectively, while the total number of binary variables in the full problem ranged from 68 for the smallest problem to 52,156 for the largest problem we were able to solve optimally. We also list the CPU times for the CPLEX solver to find the optimum solution and for our algorithm to converge, along with the percentage gap between the cost of the solution from the algorithm and the true optimum cost.

There are several conclusions that we can draw. First, the effort required to solve a problem optimally depends on the combination of factors listed in the tables (nodes, hubs, binary variables, and population density) and there is no clear direct relationship with any one specific factor. However, as might be expected, the total number of 0/1 variables seems significant. The smaller problems with under (say) 15,000 binary variables are directly solvable in a matter of seconds. However, there is one problem (no. 18) that took almost 10 hours to solve! With larger problems, these times start to go up. In at least two instances (nos. 23 and 24) the solution time was in the order of days. Yet, there was a problem (no. 6) that was larger than either of these and that could be solved in a little over 6 minutes. In general, it is hard to pinpoint what specific characteristics make the problems harder to solve optimally, and if we are presented with a new problem it would be hard to say how CPLEX might perform.

On the other hand, the iterative cyclic approach of our algorithm appears to be much more stable in its performance. Convergence is achieved in under a maximum of around two minutes for virtually all the problems tested, and the longest it took (no. 19) was about 3.5 minutes. More importantly, the solution that it finds has a cost that is always within 1% (and more often within about 0.5%) of the true minimum cost. Thus the algorithm generates high quality solutions with a substantially smaller amount of computation time than direct solution, and due to the ease with which it can be implemented, it can serve as a simple alternative to solving MIP-1 directly when a good solution is required quickly.

Finally, it is worth mentioning that we also generated much larger problems (including problems that represented the complete network for each country) but these could not be solved optimally; our algorithm also failed to converge because solving even the restricted problems om steps 2 and 3 starts to become impossible. Clearly, future research calls for other heuristic approaches, and the authors are currently working in developing these (Yang, et al., 2019).

## 6 Acknowledgment

This work was partially supported by the National Science Foundation via Award No. CMII-1536430.

## 7 References

Bland J, Clements J (1997). Protecting the world's children: the story of WHO's immunization programme. World Health Forum 19: 162-173

Chen S-I, Norman BA, Rajgopal J et al (2014). A planning model for the WHO-EPI vaccine distribution network in developing countries, IIE Trans 47:1-13

Gavi (2017). Gavi Progress Reports, Global Vaccine Alliance. https://www.gavi.org/progress-report/. Cited May 14, 2019

Gavi (2019). History of Gavi. http://www.gavi.org/about/mission/history/. Cited May 14, 2019

Lim J, Norman BA, Rajgopal J (2017). Redesign of vaccine distribution networks in low and middle-income countries, Tech Rep No. TR17-1., Department of Industrial Engineering, University of Pittsburgh, Pittsburgh.

Melkote S, Daskin MS (2001). Capacitated facility location/network design problems. Eur J Oper Res 129:481-495

Melo TM, Nickel S, Saldanha-Da-Gama F (2009). Facility location and supply chain management–A review, Eur J Oper Res 196:401-412

Şahin G, Süral H (2007). A review of hierarchical facility location models. Comp Oper Res 34:2310-2331

WHO (2013). The Expanded Programme on Immunization. http://www.who.int/immunization/programmes_systems/supply_chain/benefits_of_immunization/en/. Cited May 14, 2019

WHO (2018a). Immunization. http://www.who.int/topics/immunization/en/. Cited May 14, 2019

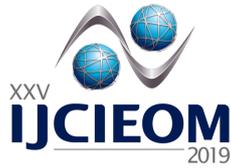

WHO (2018b). Immunization coverage. http://www.who.int/news-room/fact-sheets/detail/immunization-coverage. Cited May 14, 2019

WHO (2018c). Prequalified Devices and Equipment. http://apps.who.int/immunization_standards/vaccine_quality/pqs_catalogue/categorylist.aspx?cat_type=device. Cited May 14, 2019

Yadav P, Lydon P, Oswald J et al (2014). Integration of vaccine supply chains with other health commodity supply chains: A framework for decision making. Vaccine 32:6725-6732

Yang Y, Bidkhori H, Rajgopal J (2019). Optimizing vaccine distribution networks in low and middle-income countries, Working Paper, Department of Industrial Engineering, University of Pittsburgh, Pittsburgh.